\documentclass[12pt,a4]{article}

\usepackage{amsthm}
\usepackage{amsmath}
\usepackage{amsfonts}
\usepackage{amssymb}
\usepackage{fullpage}
\usepackage{stmaryrd}
\usepackage{graphicx}
\usepackage{psfrag}

\usepackage{pslatex}
\newtheorem{lemma}{Lemma}
\newtheorem{cor}[lemma]{Corollary}
\newtheorem{theorem}[lemma]{Theorem}

\theoremstyle{definition}
\newtheorem{defn}[lemma]{Definition}
\newtheorem*{rem}{Remark}

\newcommand{\Ref}[1]{(\ref{#1})}

\newcommand{\Diff}[2]{ \left( \frac{\partial #1}{\partial #2} \right)} 
\newcommand{\DDiff}[3]{ \left( \frac{\partial^{#3} #1}{\partial #2^{#3}} \right)} 
\newcommand{\Ddiff}[3]{ \frac{\partial^{#3} #1}{\partial #2^{#3}}} 

\newcommand{\nn}{\nonumber \\}

\newcommand{\spacebreak}{\begin{displaymath}
    \triangleleft \; \lhd \;
    \diamond \;
    \rhd \; \triangleright
  \end{displaymath}}

\newcommand{\nds}{\preceq_s}

\newcommand{\st}{\; | \;}

\newcommand{\seq}[1]{\llbracket \, #1 \, \rrbracket}

\newcommand{\hhp}[1]{{| #1|_{\Leftrightarrow}}}
\newcommand{\vhp}[1]{{| #1|_{\Updownarrow}}}

\newcommand{\hp}[1]{ \, {\odot}_{\, #1} \; }

\begin{document}

\author{Andrew Rechnitzer  \\ 
  {\small Department of Mathematics and Statistics,} \\
  {\small The University of Melbourne, Parkville Victoria 3010, Australia.}\\
  {\small email: \texttt{A.Rechnitzer@ms.unimelb.edu.au}} \\
}
\title{Haruspicy 2: \\ The anisotropic generating function of \\
self-avoiding polygons is not D-finite}
\maketitle


\begin{abstract}
  We prove that the anisotropic generating function of self-avoiding
  polygons is not a D-finite function --- proving a conjecture of
  Guttmann and Enting \cite{Guttmann2000, Guttmann1996}. This result
  is also generalised to self-avoiding polygons on hypercubic
  lattices.  Using the haruspicy techniques developed in an earlier
  paper \cite{ADR_haru}, we are also able to prove the form of the
  coefficients of the anisotropic generating function, which was first
  conjectured in \cite{Guttmann1996}.
\end{abstract}


\section{Introduction}

\begin{figure}[h!]
  \centering
  \includegraphics[scale=0.5]{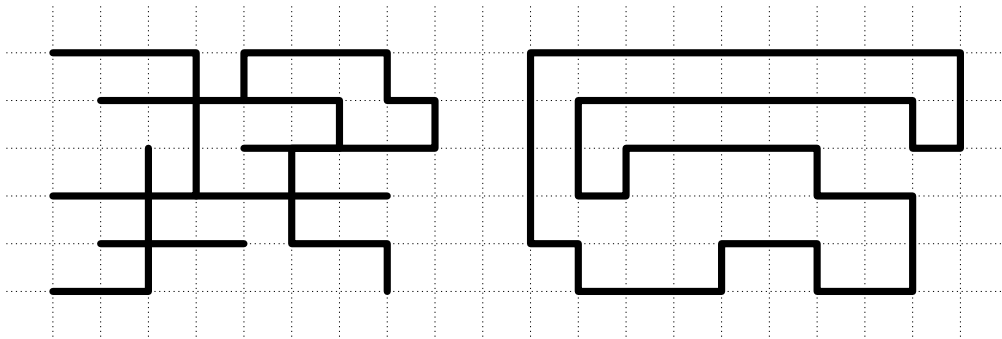}
  \caption{A square lattice bond animal (left) and a self-avoiding polygon (right).}
  \label{fig animal eg}
\end{figure}

Lattice models of magnets and polymers in statistical physics lead
naturally to questions about the combinatorial objects that form their
basis --- \emph{lattice animals}. Despite intensive study these
objects, and the lattice models from which they arise, have
tenaciously resisted rigorous analysis and much of what we know is the
result of numerical studies and ``not entirely rigorous'' results from
conformal field theory.

Recently, Guttmann and Enting \cite{Guttmann2000, Guttmann1996}
suggested a numerical procedure for testing the ``solvability'' of
lattice models based on the study of the singularities of their
\emph{anisotropic generating functions}. The application of this test
provides compelling evidence that the solutions of many of these
models do not lie inside the class of functions that includes the most
common functions of mathematical physics, namely
\emph{differentiably-finite} or \emph{D-finite} functions (defined
below). The main result of this paper is to sharpen this numerical
evidence into proof for a particular model --- \emph{self-avoiding
  polygons}.

Let us now define some of the terms we have used above. A \emph{bond
  animal} is a connected union of edges, or \emph{bonds}, on the
square lattice. The set, $\mathcal{P}$, of square lattice \emph{self-avoiding polygons},
or SAPs, is the set of all bond animals in which every vertex has
degree $2$. Equivalently it is the set of  all bond animals that are the embeddings of
a simple closed loop into the square lattice (see Figure~\ref{fig
  animal eg}).  Self-avoiding polygons were introduced in 1956 by Temperley
\cite{tempa} in work on lattice models of the phase transitions of
magnets and polymers. Not only is this problem of considerable
interest in statistical mechanics, but is an interesting combinatorial
problem in its own right. See \cite{HughesBook1, Madras93} for
reviews of this topic. 

While the model was introduced nearly 60 years ago, little progress
has been made towards either an explicit, or useful implicit,
solution. To date, only subclasses of polygons have been solved and
all of these have quite strong convexity conditions which render the
problem tractable (see \cite{MBM1996, tempa} for example).

 We wish to enumerate SAPs according to the number of
bonds they contain; since this number is always even it is customary
to count their \emph{half-perimeter} which is half the number of
bonds. The generating function of these objects is then
\begin{equation}
  P(x) = \sum_{P \in \mathcal{P}} x^{|P|},
\end{equation}
where $|P|$ denotes the half-perimeter of the polygon $P$.

To form the the \emph{anisotropic} generating function we distinguish between
vertical and horizontal bonds, and so count according to the vertical
and horizontal half-perimeters. This generating function is then
\begin{equation}
  P(x,y) = \sum_{P \in \mathcal{P}} x^\hhp{P} y^\vhp{P},
\end{equation}
where $\hhp{P}$ and $\vhp{P}$ respectively denote the horizontal and
vertical half-perimeters of $P$. By partitioning $\mathcal{P}$
according to the vertical half-perimeter we may rewrite the above
generating function as
\begin{equation}
  P(x,y) = \sum_{n \geq 1} y^n \sum_{P \in \mathcal{P}_n} x^\hhp{P} =
  \sum_{n\geq1} H_n(x) y^n,
\end{equation}
where $\mathcal{P}_n$ is the set of SAPs with $2n$ vertical bonds, and
$H_n(x)$ is its horizontal half-perimeter generating function.

The anisotropic generating function is arguably a more manageable
object than the isotropic. By splitting the set of polygons into
separate simpler subsets, $\mathcal{P}_n$, we obtain smaller pieces
which is easier to study than the whole.  If one seeks to compute or
even just understand the \emph{isotropic} generating function then one
must somehow examine \emph{all} possible configurations that can occur
in $\mathcal{P}$; this is perhaps the reason that only families of
bond animals with severe topological restrictions have been solved
(such as column-convex polygons). On the other hand, if we examine the
generating function of $\mathcal{P}_n$, then the number of different
configurations that can occur is always finite. For example, if $n=1$
all configurations are rectangles, for $n=2$ all configurations
are vertically \emph{and} horizontally convex and for  $n=3$ all
configurations are vertically \emph{or} horizontally convex.  The
anisotropy allows one to study the effect that these configurations
have on the generating function in a more controlled manner.

In a similar way, the anisotropic generating function is broken up
into separate simpler pieces, $H_n(x)$, that can be calculated exactly
for small $n$. By studying the properties of these coefficients,
rather than the whole (possibly unknown) isotropic generating
function, we can obtain some idea of the properties of the generating
function as a whole.

In many cases generating functions (and other formal power series)
satisfy simple linear differential equations; an important subclass of
such series are \emph{differentiably finite} power series; a formal
power series in $n$ variables, $F(x_1,\dots, x_n)$ with complex
coefficients is said to be \emph{differentiably finite} if for each
variable $x_i$ there exists a non-trivial differential equation:
\begin{equation}
  \label{eqn dfinite defn1}
  P_d(\mathbf{x}) \Ddiff{}{x_i}{d} F(\mathbf{x}) + \dots
  P_1(\mathbf{x}) \Ddiff{}{x_i}{} F(\mathbf{x}) + \dots
  + P_0(\mathbf{x}) F(\mathbf{x}) = 0,
\end{equation}
with $P_j$ a polynomial in $(x_1, \dots, x_n)$ with complex
coefficients \cite{Lipshitz1989}.

While no solution is known for $P(x,y)$, and certainly no equation of
the form of equation~\Ref{eqn dfinite defn1},  the first few coefficients of
$y$ may expanded numerically \cite{Tony_Iwan2000} and the following
properties were observed (up to the coefficient of $y^{14}$):
\begin{itemize}
\item $H_n(x)$ is a rational function of $x$,
\item the degree of the numerator of $H_n(x)$ is equal to the degree of its denominator.
\item the first few denominators of $H_n(x)$ (we denote them $D_n(x)$) are:
  \begin{eqnarray*}
    D_1(x) & = & (1-x) \\
    D_2(x) & = & (1-x)^3 \\
    D_3(x) & = & (1-x)^5 \\
    D_4(x) & = & (1-x)^7 \\
    D_5(x) & = & (1-x)^9 (1+x)^2 \\
    D_6(x) & = & (1-x)^{11} (1+x)^4 \\
    D_7(x) & = & (1-x)^{13} (1+x)^6 (1+x+x^2) \\
    D_8(x) & = & (1-x)^{15} (1+x)^8 (1+x+x^2)^3 \\
    D_9(x) & = & (1-x)^{17} (1+x)^{10} (1+x+x^2)^5 \\
    D_{10}(x) & = & (1-x)^{19} (1+x)^{12} (1+x+x^2)^7(1+x^2).
  \end{eqnarray*}
\end{itemize}  
Similar observations have been made for a large number of solved and
unsolved lattice models \cite{Guttmann2000} and it was noted that for
\emph{solved} models the denominators appear to only contain a small
and fixed number of different factors, while for \emph{unsolved}
models the number of different factors appears to increase with $n$.
Guttmann and Enting suggested that this pattern of increasing number
of denominator factors was the hallmark of an unsolvable problem, and
that it could be used as a test of \emph{solvability}.

In \cite{ADR_haru} we developed techniques to prove these observations
for many families of bond animals. In particular, for families of
animals that are \emph{dense} (a term we will define in the next
Section), we have the following Theorem (slightly restated for SAPs):
\begin{theorem}[from \cite{ADR_haru}]
  \label{thm Hn nature}
  If $G(x,y) = \sum_{n\geq0} H_n(x) y^n$ is the anisotropic generating function of some
  dense family of polygons, $\cal P$, then
  \begin{itemize}
  \item $H_n(x)$ is a rational function, 
  \item the degree of the numerator of $H_n(x)$ cannot be greater than the degree of its
    denominator, and
  \item the denominator of $H_n(x)$ is a product of cyclotomic polynomials.
  \end{itemize}
\end{theorem}
\begin{rem}
  We remind the reader that the \emph{cyclotomic polynomials},
  $\Psi_k(x)$, are the factors of the polynomials $(1-x^n)$. More
  precisely $(1-x^n) = \prod_{k|n} \Psi_k(x)$.
\end{rem}

In Section~\ref{sec denom bound}, we quickly review the
\emph{haruspicy} techniques developed in \cite{ADR_haru} and use them
to find a multiplicative upper bound $B_n(x)$ on the denominator
$D_n(x)$. \emph{Ie} we find a sequence of polynomials $B_n(x)$ such
that they are always divisible by $D_n(x)$.

In Section~\ref{sec dfinite}, we further refine this result to prove
that that $D_{3k-2}(x)$ contains exactly one factor of $\Psi_k(x)$
(for $k \neq 2$). This implies that the singularities of the functions
$H_n(x)$ form a dense set in the complex plane. Consequently, the
generating function $P(x,y)$ is \emph{not} differentiably finite ---
as predicted by the Guttmann and Enting solvability test. This result
is then extended to self-avoiding polygons on hypercubic lattices.



\section{Denominator Bounds}
\label{sec denom bound}
\subsection{Haruspicy}
The techniques developed in \cite{ADR_haru} allow us to determine
properties of the coefficients, $H_n(x)$, whether or not they are
known in some nice form. The basic idea is to reduce or squash the set
of animals down onto some sort of minimal set, and then various
properties of the coefficients may be inferred by examining the bond
configurations of those minimal animals. We refer to this approach as
\emph{haruspicy}; the word refers to techniques of divination based on
the examination of the forms and shapes of the organs of
animals.
\begin{figure}[h!]
  \begin{center}
    \includegraphics[scale=0.6]{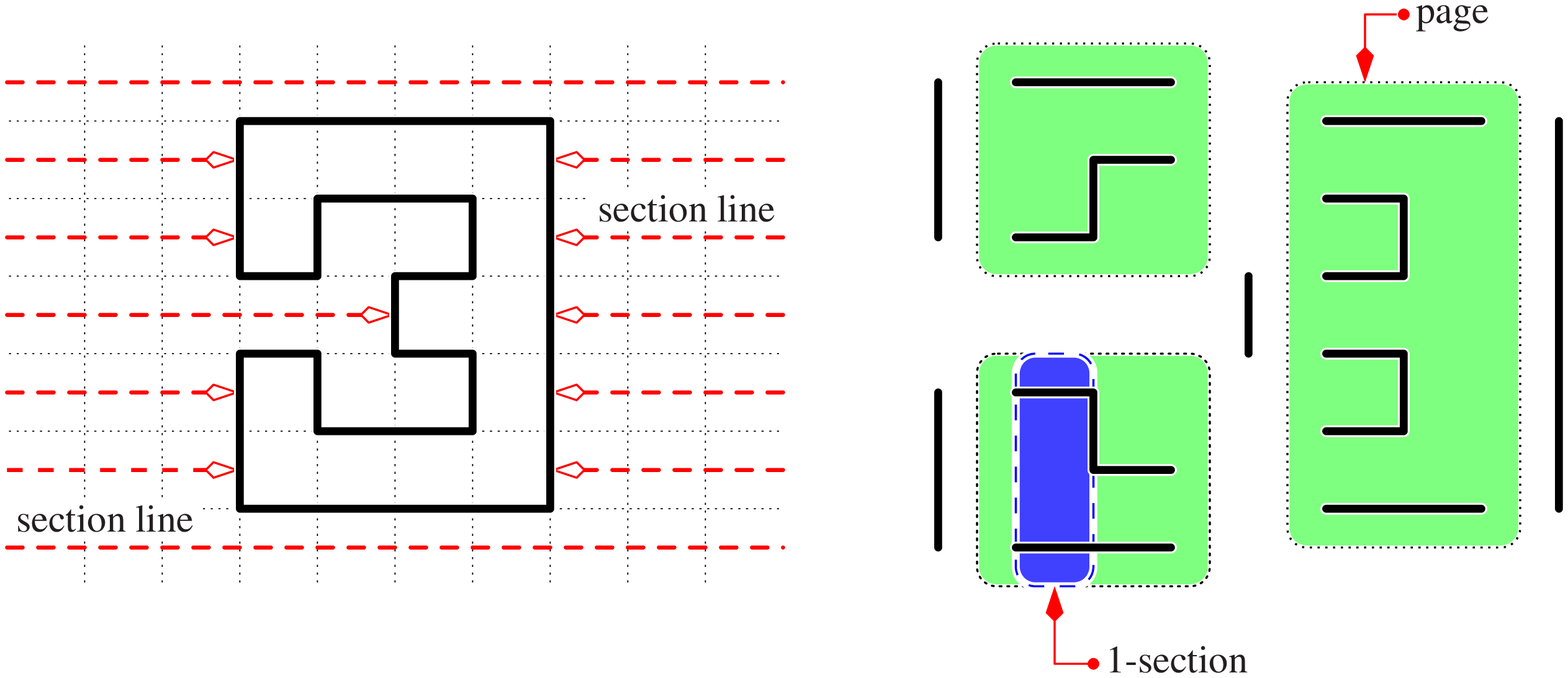}\\
    \caption{\emph{Section lines} (the heavy dashed lines in the
      left-hand figure) split the polygon into \emph{pages} (as shown
      on the right-hand figure).  Each column in a page is a
      \emph{section}. This polygon is split into 3 pages, each
      containing 2 sections; a $1$-section is highlighted.  10
      vertical bonds lie between pages and 4 vertical bonds lie within
      the pages.}
    \label{fig section defn}
  \end{center}
\end{figure}

We start by describing how to cut up polygons so that they may be
reduced or squashed in a consistent way.
\begin{defn}
  Draw horizontal lines from the extreme left and the extreme right of
  the lattice towards the animal so that the lines run through the
  middle of each lattice cell.  These lines are called \emph{section
    lines}.  The lines are terminated when they first touch (\emph{ie}
  are obstructed by) a vertical bond (see Figure~\ref{fig section
    defn}).
  
  Cut the lattice along each section line from infinity until it
  terminates at a vertical bond. Then from this vertical bond cut
  vertically in both directions until another section line is reached.
  In this way the polygon (and the lattice) is split into \emph{pages}
  (see Figure~\ref{fig section defn}); we consider the vertical bonds
  along these vertical cuts to lie \emph{between} pages, while the
  other vertical bonds lie \emph{within} the pages.
  
  We call a \emph{section} the set of horizontal bonds within a single
  column of a given page.  Equivalently, it is the set of horizontal
  bonds of a column of an animal between two neighbouring section
  lines. A section with $2k$ horizontal bonds is a $k$-section.  The
  number of $k$-sections in a polygon, $P$, is denoted by
  $\sigma_k(P)$.
\end{defn}

\begin{figure}[ht]
  \begin{center}
    \includegraphics[scale=0.55]{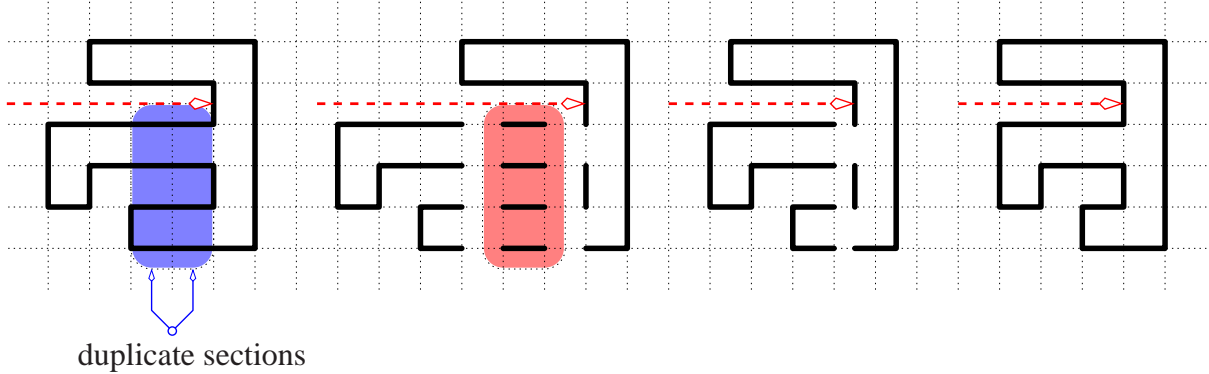}
    \caption{The process of section deletion. The two indicated sections 
      are identical. Slice either side of the duplicate and separate the polygon into
      three pieces. The middle piece, being the duplicate, is removed and the remainder of
      the polygon is recombined. Reversing the steps leads to section
      duplication. Also indicated is a section line which separates
      the duplicate sections from the rest of the columns in which
      they lie.}
    \label{fig explicit section delete}
  \end{center}
\end{figure}

The polygon has now been divided up into smaller units, which we have
called sections. In some sense many of these sections are superfluous
and are not needed to encode its ``\emph{shape}'' (in some loose sense
of the word). More specifically, if there are two neighbouring
sections that are the same, then we can reduce the polygon by removing
one of them, while leaving the polygon with essentially the same shape.
\begin{defn}
  We say that a section is a \emph{duplicate section} if the section
  immediately on its left (without loss of generality) is identical
  (see Figure~\ref{fig explicit section delete}).
  
  One can reduce polygons by \emph{deletion} of duplicate sections by
  slicing the polygon on either side of the duplicate section,
  removing it and then recombining the polygon, as illustrated in
  Figure~\ref{fig explicit section delete}.  By reversing the
  section-deletion process we define \emph{duplication} of a section.

  We say that a set of polygons, $\mathcal{P}$, is \emph{dense} if the
  set is closed under section deletion and duplication. \emph{Ie} no
  polygon outside the set can be produced by section deletion and / or
  duplication from a polygon inside the set.
\end{defn}

The process of section-deletion and duplication leads to a partial
order on the set of polygons.
\begin{defn}
  For any two polygons $P, Q \in {\cal P}_n$, we write $P \nds Q$ if
  $P=Q$ or $P$ can be obtained from $Q$ by a sequence of
  section-deletions.  A \emph{section-minimal} polygon, $P$, is a
  polygon such that for all polygons $Q$ with $Q \nds P$ we have
  $P=Q$.
\end{defn}
The lemma below follows from the above definition (see \cite{ADR_haru} for
details):
\begin{lemma}
  The binary relation $\nds$ is a partial order on the set of
  polygons. Further every polygon reduces to a unique section-minimal
  polygon, and there are only a finite number of minimal polygons in
  $\mathcal{P}_n$.
\end{lemma}

By considering the generating function of all polygons that are equivalent (by some
sequence of section-deletions) to a given section-minimal polygon, we
find that $H_n(x)$ may be written as the sum of simple rational
functions. Theorem~\ref{thm Hn nature} follows directly from
this. Further examination of the denominators of these functions gives
the following result:
\begin{theorem}[from \cite{ADR_haru}]
  \label{thm poles sections}
  If $H_n(x)$ has a denominator factor $\Psi_k(x)$, then ${\cal P}_n$ must
  contain a section-minimal polygon containing a $K$-section for some $K \in
  \mathbb{Z}^+$ divisible by $k$. Further if $H_n(x)$ has a denominator factor
  $\Psi_k(x)^{\alpha}$, then ${\cal P}_n$ must contain a section-minimal polygon
  that contains $\alpha$ sections that are $K$-sections for some (possibly
  different) $K \in \mathbb{Z}^+$ divisible by $k$.
\end{theorem}
This theorem demonstrates the link between the factors of $D_n(x)$ and
the sections in section-minimal polygons with $2n$ vertical bonds.

\subsection{The number of $k$-sections}
In this subsection, we shall demonstrate the following multiplicative
upper bound on the denominator, $D_n(x)$ of $H_n(x)$:
\begin{eqnarray}
  D_n(x) & \mbox{ is a factor of } & \prod_{k=1}^{\lceil n/3 \rceil}
    \Psi_k(x)^{2n-6k+5} .
\end{eqnarray}
We do this by finding an upper bound on the number of $k$-sections
that a SAP with $2n$ vertical bonds may contain. A proof of the
corresponding result for general bond animals is given in
\cite{ADR_haru}; here we follow a similar line of proof, but
specialise (where possible) to the case of SAPs.

The proof consists of several steps:
\begin{enumerate}
\item Find the maximum number of sections in a polygon with $2V$
  vertical bonds.
\item Determine a lower bound on the number of vertical bonds and
  sections that must lie to the left (without loss of generality) of a
  $k$-section. This gives a lower bound on the number of sections that
  must lie to the left of the leftmost $k$-section and to the right of
  the rightmost $k$-section --- none of these can be $k$-sections and
  so we obtain a lower bound on number of sections that cannot be
  $k$-sections.
\item From the above two facts we obtain an upper bound on the number
  of sections in a polygon that may be $k$-sections; assume that they
  are all $k$-sections.
\item Theorem~\ref{thm poles sections} then gives the upper bound on
  the exponent of $\Psi_k(x)$.
\end{enumerate}

Please note that for the remainder of this part of the paper we shall
use ``\emph{sm-polygons}'' to denote ``\emph{section-minimal polygons}''
unless otherwise stated.

\begin{lemma}
  \label{lem num sec in page}
  An sm-polygon that contains $p$ pages and $v$ vertical bonds inside
  those pages may contain at most $p+v$ sections.
\end{lemma}
\proof Consider the $v_i$ vertical bonds inside the $i^{th}$
page. Between any two sections in this page there must be at least 1
vertical bond (otherwise the horizontal bonds in both sections would
be the same and they would be duplicate sections). Hence the $i^{th}$
page contains at most $v_i+1$ sections. Since every section must lie
in exactly 1 page the result follows. \qed

\begin{lemma}
  \label{lem pages and rows}
  The maximum number of pages in an sm-polygon is $2R-1$
  where $R$ denotes the number of rows in the sm-polygon.
\end{lemma}
\proof See Lemma~13 in \cite{ADR_haru}. We note that this differs from the result
for bond animal since all sections must contain an even number of
horizontal bonds, and must also lie between vertical bonds.
Consequently we are only interested in those pages that lie
\emph{inside} the sm-polygon. \qed

\begin{lemma}
  The maximum number of sections in an sm-polygon with $2V$ vertical bonds
  is $2V-1$.
  \label{lem max sec}
\end{lemma}
\proof Consider an sm-polygon of height $R$ with $2V = 2R+2v$ vertical
bonds.  Of these vertical bonds $2R$ block section lines and the
remaining $2v$ may lie \emph{inside} pages. By Lemma~\ref{lem pages
  and rows}, this sm-polygon has at most $2R-1$ pages. At most $2v$
vertical bonds lie inside these pages and so by Lemma~\ref{lem num sec
  in page} the result follows.  \qed

\vspace{0.5cm}

The above lemma tells us the maximum number of sections in a
sm-polygon. We now determine how many of these sections lie to the
left (without loss of generality) of a $k$-section. We start by
determining how many vertical bonds lie to the left of a $k$-section.

\begin{figure}[ht]
  \begin{center}
    \includegraphics[scale=0.55]{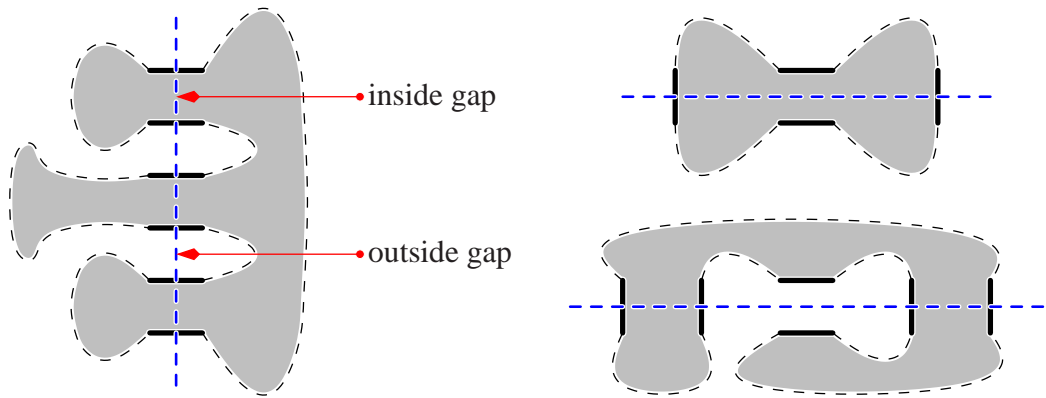}
    \caption{Vertical and horizontal lines drawn through a $k$-section
      show the minimum number of vertical bonds required in their
      construction.}
    \label{fig cut ksect}
  \end{center}
\end{figure}

\begin{lemma}
  To the left (without loss of generality) of a $k$-section there are
  at least $3k-2$ vertical bonds, of which at least $2k-1$ obstruct
  section lines.
  
  Hence no polygon with fewer than $6k-4$ vertical bonds may contain a
  $k$-section. Further, it is always possible to construct a polygon
  with $6k-4$ vertical bonds and a single $k$-section.
  \label{lem first ksec}
\end{lemma}
\proof Consider a vertical line drawn through a $k$-section (as
depicted in left-hand side of Figure~\ref{fig cut ksect}). The line
starts outside the polygon and then as it crosses horizontal bonds it
alternates between the inside and outside of the polygon.  More
precisely, there are $k+1$ segments of the line that lie outside the
polygon and $k$ segments that lie inside the polygon. Let us call the
segments that lie within the polygon ``\emph{inside gaps}'' and those
that lie outside ``\emph{outside gaps}''.

Draw a horizontal line through an inside gap (as depicted in the
top-right of Figure~\ref{fig cut ksect}). If we traverse the
horizontal line from left to right, we must cross at least 1 vertical
bond to the left of the gap (since it is inside the polygon) and then
another to the right of the gap. Hence to the left of any inside gap
there must be at least 1 vertical bond. Similarly we must cross at
least 1 vertical bond to the right of any inside gap.

Draw a horizontal line through the topmost of the $k+1$ outside gaps.
Since the line need not intersect the polygon it need not cross any
vertical bonds at all. Similarly for the bottommost outside gap.

Now draw a horizontal line through one of the other outside gaps (as
depicted in the bottom-right of Figure~\ref{fig cut ksect}). Traverse
this line from the left towards the outside gap. If no vertical bonds
are crossed then a section line may be drawn from the left into the
outside gap. This splits the $k$-section into two smaller sections
and so contradicts our assumptions. Hence we must cross at least 1
vertical bond to block section lines. If we cross only a single
vertical bond before reaching the gap then the gap would lie inside
the polygon. Hence we must cross at least 2 (or any even number)
vertical bonds before reaching the gap.  Similar reasoning shows that
we must also cross an even number of vertical bonds when we continue
traversing to the right.

Since any $k$-section contains $k$ inside gaps, a topmost outside gap,
a bottommost outside gap and $k-1$ other outside gaps, there must be
at least $k\times1+2\times0 + 2\times(k-1) = 3k-2$ vertical bonds to its
left and $3k-2$ vertical bonds to its right.

\begin{figure}[ht]
  \begin{center}
    \includegraphics[scale=0.44]{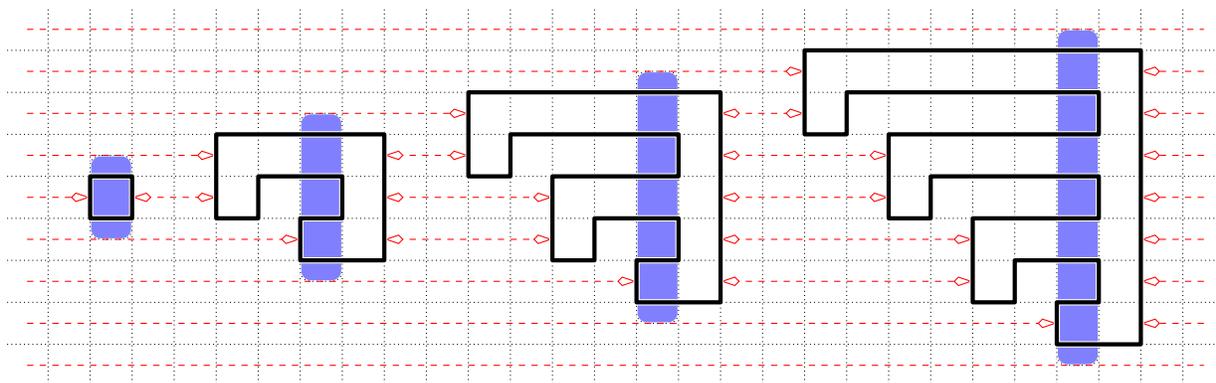}
    \caption{Section-minimal polygons with $6k-4$
      vertical bonds and a single $k$-section may be constructed by
      concatenating such ``\emph{hook}'' configurations.}
    \label{fig min ksec}
  \end{center}
\end{figure}

Consider the polygons depicted in Figure~\ref{fig min ksec} that are
constructed by adding ``\emph{hooks}''. In this way we are always
able to construct an sm-polygon with $(6k-4)$ vertical bonds and
exactly one $k$-section. \qed

\begin{figure}[ht]
  \centering

  \psfrag{1}[Bl][][0.7]{1}  \psfrag{2}[Bl][][0.7]{2}  \psfrag{3}[Bl][][0.7]{3}
  \psfrag{4}[Bl][][0.7]{4}  \psfrag{5}[Bl][][0.7]{5}  \psfrag{6}[Bl][][0.7]{6}
  \psfrag{7}[Bl][][0.7]{7}  \psfrag{8}[Bl][][0.7]{8}  \psfrag{9}[Bl][][0.7]{9} 
  \includegraphics[scale=0.33]{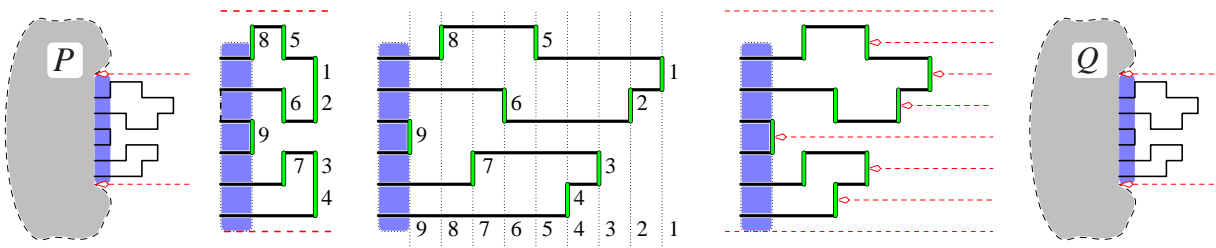}
  \caption{Given an sm-polygon $P$ we can create a new
    sm-polygon $Q$ that is identical to $P$ except for the region
    lying to the right of its rightmost $k$-section; that region is
    altered to maximise the number of sections lying to the right of
    the $k$-section. We do this by stretching that portion $P$ that
    lies to the right of the $k$-section so that no two vertical bonds
    lie in the same vertical line. The polygon is then made
    section-minimal again by deleting duplicate sections.}
  \label{fig stretch}
\end{figure}

The next lemma shows that given an sm-polygon, $P$, that contains a
$k$-section, we are always able to find a new sm-polygon, $Q$,
with the same number of vertical bonds that has \emph{at least}
$3(k-1)$ sections to the left of its leftmost $k$-section. This result
allows us to compute how many sections in an sm-polygon cannot be
$k$-sections since they lie to the left of the leftmost or to the
right of the rightmost $k$-section.

\begin{lemma}
  Let $P$ be an sm-polygon that contains a $k$-section and $2V$
  vertical bonds. If there are fewer than $3(k-1)$ sections to the
  right of the rightmost $k$-section in $P$, then there exists another
  sm-polygon, $Q$, that is identical to $P$ except that to the right
  of the rightmost $k$-section there are at least $3(k-1)$
  sections. See Figure~\ref{fig stretch}.
  
  Similarly given a polygon, $P'$ with fewer than $3(k-1)$ sections to
  the left of the leftmost $k$-section, there exists another polygon
  $Q'$ identical to $P'$ except that there are at least $3(k-1)$
  sections to the left of the leftmost $k$-section.
  \label{lem stretch}
\end{lemma}
\proof We prove the above result by ``\emph{stretching}'' the portion
of the sm-polygon, $P$, to the right of the rightmost $k$-section so
as to obtain a new sm-polygon, $Q$, in which the number of sections to
the right of the $k$-section is maximised.

Consider the example given in Figure~\ref{fig stretch}. Consider the
portion of the sm-polygon that lies to the right of rightmost $k$-section
(which is highlighted). Label the vertical bonds from top-rightmost
(1) to bottom-leftmost (9). We now ``\emph{stretch}'' the horizontal
bonds of the sm-polygon so that bonds with higher labels lie to the left
of those with lower labels and so that no two bonds lie in the same
vertical line (Figure~\ref{fig stretch}, centre). To recover a
section-minimal polygon we now delete duplicate sections
(Figure~\ref{fig stretch}, right). We now need to determine how many
sections remain.

\begin{figure}[h]
  \centering
  \includegraphics[scale=0.45]{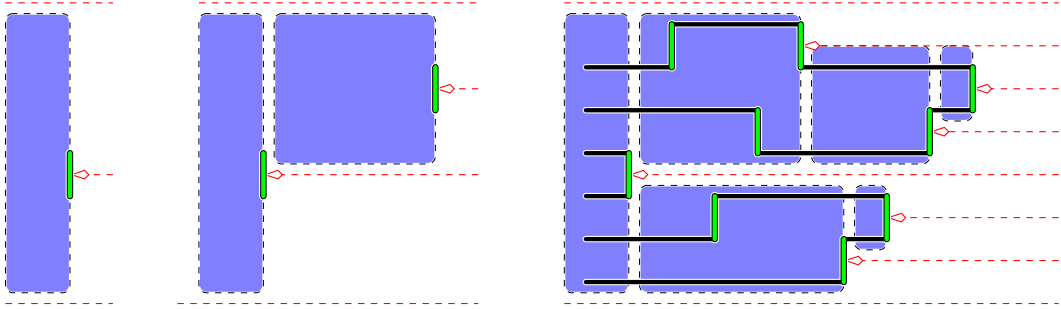}
  \caption{The pages in the stretched polygon before removing
    duplicate sections. By ``\emph{scanning}'' from left to right we see that
    each page corresponds one vertical bond that blocks a
    section line.}
  \label{fig stretch 2}
\end{figure}

Consider the stretched portion of polygon before duplicate sections
are removed.  If there were originally $r$ vertical bonds blocking
section lines, then there are still $r$ vertical bonds blocking
section lines after stretching. See Figure~\ref{fig stretch 2}. Since
no two vertical bonds lie in the same vertical line, each page
corresponds to a single vertical bond that blocks a section line
(which will lie on the right-hand edge of the page).  Hence the
stretched portion polygon contains $r$ pages (one of which contains
the $k$-section). The other vertical bonds must lie within these
pages. See also the proof of Lemma~13 in \cite{ADR_haru}.

Thus, if there were $v=r+m$ vertical bonds to the right of the
$k$-section, with $r$ blocking section lines, then after deleting
duplicate sections there will be $r$ pages (no pages will be removed)
and $m$ vertical bonds within those pages (with no two vertical bonds
in the same page lying in the same vertical line). Consequently there
will be $r+m-1$ sections \emph{excluding} the $k$-section.

By Lemma~\ref{lem first ksec} there must be at least $3k-2$ vertical
bonds to the right of a $k$-section , and so the ``stretching''
procedure will produce an sm-polygon with at least $3k-3$ sections to the
right of the $k$-section.  

Note that this procedure does not change the number of vertical bonds
in each row, nor the number of vertical bonds on either side of the
$k$-section.  \qed

Since we now know the total number of sections in a section minimal
polygon and how many of these cannot be $k$-sections we can prove an
upper bound on the number of $k$-sections:
\begin{theorem}
  A section-minimal polygon $P$ that contains $2V = (6k-4+2M)$
  vertical bonds may not contain more than $2M+1$ sections that
  contain $2k$ or more horizontal bonds.
  \label{thm num ksec}
\end{theorem}
\proof By Lemma~\ref{lem max sec}, $P$ may contain no more that
$(2V-1)$ sections. We complete the proof by assuming that the theorem
is false and then reaching a contradiction.

Consider an sm-polygon, $P$, that does not have a section with $>2k$
horizontal bonds, but does contain more than $(2M+1)$ $k$-sections. By
Lemma~\ref{lem stretch} we may always ``\emph{stretch}'' the portion
of the polygon lying to the right of the rightmost $k$-section to obtain
a new section-minimal polygon so that at least $3(k-1)$ sections lie
to the right of the rightmost $k$-section. Similarly we may ``stretch''
the portion of the polygon lying to the left of the leftmost
$k$-section to obtain a new section-minimal polygon $Q$ that has at
least $6(k-1)$ sections lying either to the left of the leftmost or to
the right of the rightmost $k$-sections. Consequently this new polygon
contains more than $(2M+1)+(6k-6) = 6k+2M-5$ sections. This
contradicts Lemma~\ref{lem max sec}.

Now consider an sm-polygon that contains sections with at least $2k$
horizontal bonds. Assume that it does contain more than $2M+1$ such
sections. Without loss of generality consider the leftmost section
with at least $2k$ horizontal bonds.  By applying Lemma~\ref{lem
  stretch} we see that one may always construct a new section-minimal
polygon so that at least $3(k-1)$ sections lie to the left of the
leftmost such section.  Repeating the argument in the paragraph above
shows that one will reach a contradiction and the proof is complete.
\qed

\begin{figure}[h!]
  \begin{center}
    \includegraphics[scale=0.5]{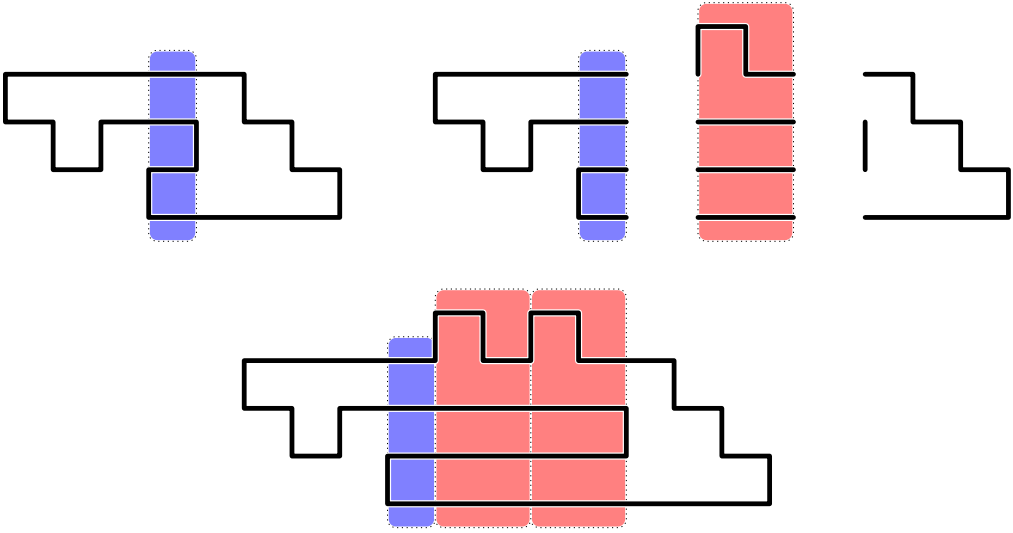}
    \caption{To construct a  polygon with  $2M+1$ $k$-sections and $6k-4+2M$ vertical
      bonds, start with a polygon with a single $k$-section and $6k-4$
      vertical bonds as shown (top left). Cut it on the right of the
      $k$-section. Insert $M$ copies of the pair of $k$-sections and
      recombine the polygon.  This gives a polygon with $2M+1$
      $k$-sections and $6k-4+2M$ vertical bonds.}
    \label{fig most ksec}
  \end{center}
\end{figure}

\begin{rem}
  It is possible to construct a section-minimal polygon with exactly
  $(6k-4+2M)$ vertical bonds and $2M+1$ $k$-sections --- see
  Figure~\ref{fig most ksec}.
\end{rem}

\begin{cor}
  \label{cor denom bound}
  The factor of $\Psi_k(x)$ in the denominator, $D_n(x)$ of $H_n(x)$ may not
  appear with a power greater than $2n-6k+5$. Hence we have the
  following multiplicative upper bound for $D_n(x)$:
  \begin{displaymath}
    D_n(x) \; \left| \; \prod_{k=1}^{\lceil n/3 \rceil} \Psi_k(x)^{2n-6k+5} \right. .
  \end{displaymath}
\end{cor}
\proof This follows by combining the results of Theorems~\ref{thm Hn
  nature},~\ref{thm poles sections}~and~\ref{thm num ksec}.

\begin{rem}
  We note that a comparison of the above bound on the denominator of
  $H_n(x)$ appears to be quite tight when compared with series
  expansion data~\cite{Tony_Iwan2000}. It appears to be wrong only by a single
  factor of $\Psi_2(x)$; the exponents of other factors appear to be
  equal to that of the bound.
  
  We also note that the above result significantly reduces the
  difficultly of computing the coefficients, $H_n(x)$ of the
  anisotropic generating function. In particular, we know that
  $H_n(x)$ is a rational function whose numerator degree is no greater
  than that of its denominator.  Corollary~\ref{cor denom bound} gives
  this denominator (up to multiplicative cyclotomic factors) and as a
  consequence also bounds the degree of the numerator and hence the
  number of unknowns we must compute in order to know $H_n(x)$.
  
  Since the degree of $\Psi_k(x)$ is no greater than $k$, the degree
  of $D_n(x)$ is no greater than $\sum_{k=1}^{\lceil n/3 \rceil}
  k(2n-6k+5) \sim \frac{1}{27} n^3$. Note that using similar
  (non-rigorous) arguments to those in Section~4.2 of \cite{ADR_haru}
  one can show that the degree grows like $\frac{2}{9\pi^2} n^3$.
  Hence the degree of the numerator (and the number of unknowns to be
  computed) grows as $n^3$. Bounds from transfer matrix techniques
  (such as \cite{Enting1980}) grow exponentially.
\end{rem}

\section{The nature of the generating function}
\label{sec dfinite}
\subsection{Differentiably finite functions}
Perhaps the most common functions in mathematical physics (and
combinatorics) are those that satisfy simple linear differential
equations. A subset of these are the differentiably finite functions
that satisfy linear differential equations with polynomial coefficients.

\begin{defn}
  \label{defn dfinite}
  Let $F(x)$ be a formal power series in $x$ with coefficients in
  $\mathbb{C}$.  It is said to be \emph{differentiably finite} or
  \emph{D-finite} if there exists a non-trivial differential equation:
  \begin{equation}
    P_d(x) \Ddiff{}{x}{d} F(x) + \dots
    + P_1(x) \Ddiff{}{x}{} F(x) 
    + P_0(x) F(x) = 0,
  \end{equation}
  with $P_j$ a polynomial in $x$ with complex coefficients
  \cite{Lipshitz1989}. 
  
  In this paper we consider series, $G(x,y)$ that are formal power
  series in $y$ with coefficients that are rational functions of $x$. Such
  a series is said to be D-finite if there exists a non-trivial
  differential equation:
  \begin{equation}
    \label{eqn dfinite de}
    Q_d(x,y) \Ddiff{}{y}{d} G(x,y) + \dots
    + Q_1(x,y) \Ddiff{}{y}{} G(x,y) 
    + Q_0(x,y) G(x,y) = 0,
  \end{equation}
  with $Q_j$ a polynomial in $x$ and $y$ with complex coefficients
\end{defn}

One of the main aims of this paper is to demonstrate that the
anisotropic generating function of SAPs is not D-finite, and we do so
by examining the singularities of that function.

The classical theory of linear differential equations implies that a
D-finite power series of a single variable has only a finite number of
singularities. This forms a very simple ``D-finiteness test'' --- a function
such as $f(x) = 1 / \cos(x)$ cannot be D-finite since it has an
infinite number of singularities. Unfortunately we know very little
about the singularities of the \emph{isotropic} SAP generating
function and cannot apply this test.

When we turn our attention to the anisotropic generating function (a
power series with rational coefficients) there is a similar test that
examines the singularities of the \emph{coefficients}.  Consider the
following example:
\begin{equation}
  f(x,y) = \sum_{n \geq 1} \frac{x^n}{1-nx} \; y^n.
\end{equation}
The coefficient of $y^n$ is singular at $x=1/n$ and so the set of
singularities of its coefficients, {$\{ n^{-1} \st n \in \mathbb{Z}^+ \}$}, is
infinite and has an accumulation point at $0$. In spite of this it is
a D-finite power series in $y$, since it satisfies the following
partial differential equation:
\begin{equation}
  xy^2(1-xy) \Ddiff{f}{y}{2}-y(1-xy+x^2y)\Diff{f}{y}+f = 0.
\end{equation}
So the set of the singularities of the coefficients of a D-finite
series may be infinite and have accumulation points. It may not,
however, have an infinite number of accumulation points.
\begin{theorem}[from \cite{MBM_AR_CH}]
  Let $f(x,y) = \sum_{n\geq0} y^n H_n(x)$ be a D-finite series in $y$
  with coefficients $H_n(x)$ that are rational functions of $x$.
  For $n\geq0$ let $S_n$ be the set of poles of $H_n(x)$, and let $S =
  \bigcup_n S_n$. Then $S$ has only a finite number of accumulation
  points.
\label{thm Dfinite poles}
\end{theorem}

In order to apply this theorem to the self-avoiding polygon generating
function we need to prove that the denominators of the coefficients
$H_n(x)$ suggested by Corollary~\ref{cor denom bound} do not cancel
with the numerators --- so that the singularities suggested by those
denominators really do exist. Unfortunately, we are unable to prove
such a strong result. However, we do not need to understand the full
singularity structure of the coefficients; the following result is
sufficient:
\begin{theorem}
  \label{thm first ksec exp}
  For $k \neq 2$ the generating function $H_{3k-2}(x)$ has simple
  poles at the zeros of $\Psi_k(x)$. Equivalently the denominator of
  $H_{3k-2}(x)$ contains a single factor of $\Psi_k(x)$ which does not
  cancel with the numerator.
\end{theorem}
An immediate corollary of this result is that singularities of the
coefficients $H_n(x)$ are dense on the unit circle, $|x|=1$ and so the
anisotropic generating function is not a D-finite power series in $y$.

\subsection{\emph{2-4-2}~polygons}
In order to prove Theorem~\ref{thm first ksec exp} we split the set of
polygons with $(6k-4)$ vertical bonds into two sets --- polygons that
contain a $k$-section and those that do not. Let us denote those
polygons with $(6k-4)$ vertical bonds and at least one $k$-section by
$\mathcal{K}_{3k-2}$. Hence we may write the generating function
$H_{3k-3}(x)$ as
\begin{displaymath}
  H_{3k-2}(x) = \sum_{P \in \mathcal{K}_{3k-2}} x^\hhp{P} 
  + \sum_{P \in \mathcal{P}_{3k-2} \setminus \mathcal{K}_{3k-2}} x^\hhp{P}.
\end{displaymath}

\begin{lemma}
  \label{lem pole exactly one ksec}
  The factor $\Psi_k(x)$ appears in the denominator of the generating
  function\\ $\sum_{P \in \mathcal{K}_{3k-2}} x^\hhp{P}$ with exponent
  exactly equal to one if and only if it appears in the denominator
  of  $H_{3k-2}(x)$ with exponent exactly equal to one.
\end{lemma}
\proof The sets $\mathcal{K}_{3k-2}$ and $\mathcal{P}_{3k-2} \setminus
\mathcal{K}_{3k-2}$ are trivially dense, and so by Theorem~\ref{thm Hn
  nature} we know that the horizontal half-perimeter generating
functions of these sets are rational and that their denominators are
products of cyclotomic factors. Further, since $\mathcal{P}_{3k-2}
\setminus \mathcal{K}_{3k-2}$ does not contain a polygon with
$k$-section (or indeed, by Lemma~\ref{lem first ksec}, any section
with more than $2k$ horizontal bonds), it follows by Theorem~\ref{thm
  poles sections} that the denominator of the horizontal
half-perimeter generating function of this set is a product of
cyclotomic polynomials $\Psi_j(x)$ for $j$ \emph{strictly less} than
$k$. Consequently this generating function is not singular at the
zeros of $\Psi_k(x)$. By Theorem~\ref{thm num ksec}, every
section-minimal polygon in $\mathcal{K}_{3k-2}$ contains exactly one
$k$-section, and so the exponent of $\Psi_k(x)$ in in the denominator
of the horizontal half-perimeter generating function of
$\mathcal{K}_{3k-2}$ is either one or zero (due to cancellations with
the numerator). The result follows since this denominator factor may
not be cancelled by adding the other generating function. \qed

\spacebreak

The above Lemma shows that to prove Theorem~\ref{thm first ksec exp}
it is sufficient to prove a similar result for the set of polygons,
$\mathcal{K}_{3k-2}$. Let us examine this set further. In the proof of
Lemma~\ref{lem first ksec} it was shown that a $k$-section could be
decomposed an alternating sequence of ``\emph{inside gaps}'' and
``\emph{outside gaps}''; a row containing an inside gap contained at
least $2$ vertical bonds and a row containing an outside gap contained
at least $4$ vertical bonds (see the examples in Figure~\ref{fig min
  ksec}). We now concentrate on polygons containing $2$ vertical bonds
in very second row and $4$ vertical bonds in every other row.


\begin{defn}
  Number the rows of a polygon $P$ starting from the topmost row (row
  1) to the bottommost (row r). Let $v_i(P) $ be the number of
  vertical bonds in the $i^{th}$ row of $P$.  If  $(v_1(P), \dots,
  v_r(P)) \\ = (2,4,2,\dots,4,2)$ then we call $P$ a \emph{2-4-2}
  polygon.  We denote the set of such \emph{2-4-2} polygons with $2n$
  vertical bonds by $\mathcal{P}^{242}_n$. Note that this set is empty
  unless $2n=6k-4$ (for some $k = 1, 2, \dots$).

\end{defn}

\begin{lemma}
  A section-minimal polygon with $(6k-4)$ vertical bonds that contains
  one $k$-section must be a \emph{2-4-2} polygon. On the other hand, a
  section-minimal \emph{2-4-2} polygon need not contain a $k$-section.
\end{lemma}
\proof The first statement follows by arguments given in the proof of
Lemma~\ref{lem first ksec}. The rightmost polygon in Figure~\ref{fig
  242 eg} show that a \emph{2-4-2} polygon need not contain a $k$-section. \qed
\begin{figure}[h!]
  \centering
  \includegraphics[scale=0.45]{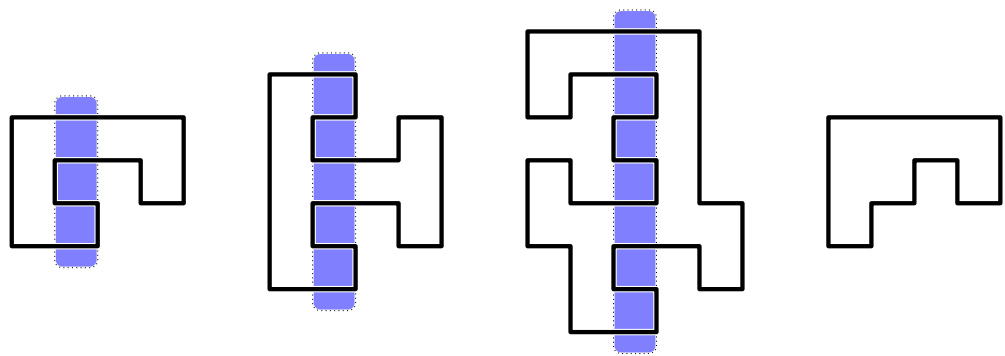}
  \caption{Four section-minimal \emph{2-4-2} polygons. The first three
    polygons contain a $2$-section, a $3$-section, and a $4$-section
    respectively. The rightmost polygon contains only $1$-sections.}
  \label{fig 242 eg}
\end{figure}

Despite the fact that \emph{2-4-2} polygons are a superset of those
polygons containing at least one $k$-section, it turns out both that
they are easier to analyse (in the work that follows) and that a
result analogous to Lemma~\ref{lem pole exactly one ksec} still holds.
\begin{lemma}
  \label{lem 242 pole}
  The factor $\Psi_k(x)$ appears in the denominator of the generating
  function $\sum_{P \in \mathcal{P}^{242}_{3k-2}} x^\hhp{P}$ with exponent
  exactly equal to $1$  if and only if it appears in the denominator
  of  $H_{3k-2}(x)$ with exponent exactly equal to one.
\end{lemma}
\proof Similar to the the proof of Lemma~\ref{lem pole exactly one
  ksec}. \qed

In the next section we derive a (non-trivial) functional equation
satisfied by the generating function of \emph{2-4-2}~polygons.

\subsection{Counting with Hadamard products}
By far the most well understood classes of square lattice polygons are
families of \emph{row convex} polygons. Each row of a row convex
polygon contains only two vertical bonds; this allows one to find a
construction by which polygons are built up \emph{row-by-row}. This
technique is sometimes called the \emph{Temperley method}
\cite{MBM1996, tempa}.

Since every second row of a \emph{2-4-2}~polygon contains $2$ vertical
bonds, we shall find a similar construction that instead of building
up the polygons row-by-row, we build them two rows at a time (an idea
also used in \cite{MBM_ADR_02}).  Like the constructions given in
\cite{MBM1996}, this construction leads quite naturally to a
functional equation satisfied by the generating function. One could
also derive this functional equation using the techniques described in
\cite{MBM1996}, however it proves more convenient in this case to use
techniques based on the application of Hadamard products (this idea is
also used in \cite{MBM1999}).

\begin{figure}[h!]
  \centering
  \includegraphics[scale=0.3]{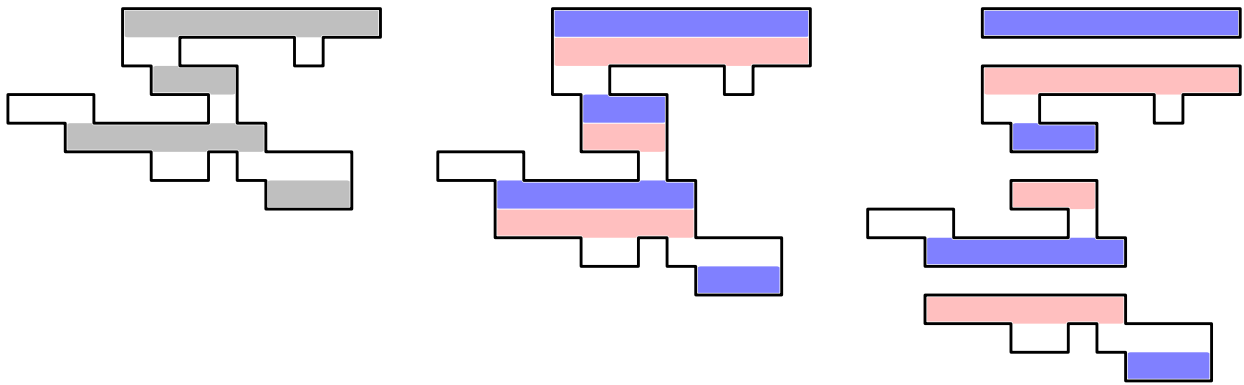}
  \caption{Decomposing \emph{2-4-2} polygons into building blocks. Highlight
    each row with $2$ vertical bonds. Then ``duplicate'' each of these
  rows excepting the bottommost. By cutting along each of these
  duplicated rows each \emph{2-4-2} polygon is decomposed into a
  rectangle (of unit height) and a sequence of building blocks.}
  \label{fig 242 decomp}
\end{figure}

We shall start by showing how \emph{2-4-2} polygons may be decomposed
into smaller units we shall call \emph{seeds} and \emph{building
  blocks}. Consider the \emph{2-4-2} polygon in Figure~\ref{fig 242
  decomp}. Start by highlighting each row with $2$ vertical bonds. We
then ``\emph{duplicate}'' each of these rows, excepting the
bottommost; this situation is depicted in the middle polygon in
Figure~\ref{fig 242 decomp}. By cutting the polygon horizontally
between each pair of duplicate rows we decompose the polygon
\emph{uniquely} into a rectangle of unit height and a sequence of
\emph{2-4-2} polygons of height $3$, such that the bottom row of each
polygon is the same length of the top row of the next in the sequence.
We refer to this initial rectangle as the \emph{seed block} and the
subsequent \emph{2-4-2} polygons of height $3$ as \emph{building
  blocks}.

\begin{figure}[ht]
  \centering
  \includegraphics[scale=0.3]{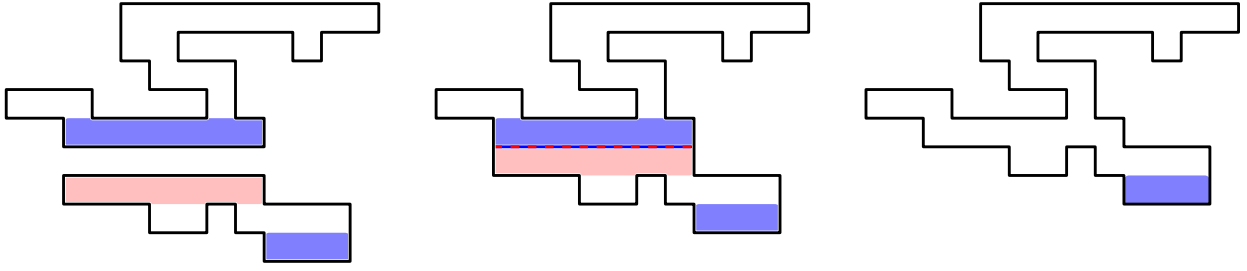}
  \caption{Constructing a \emph{2-4-2} polygon from a (shorter)
    \emph{2-4-2} polygon and a \emph{2-4-2} building block. Note that
    when the building block and the polygon are squashed together, the
    total vertical perimeter is reduced by $2$, and the total
    horizontal perimeter is reduced by twice the width of the joining
    row.}
  \label{fig 242 recomp}
\end{figure}

This decomposition implies that each \emph{2-4-2} polygon is either a
rectangle of unit height, or may be constructed by
``\emph{combining}'' a (shorter) \emph{2-4-2} polygon and a
\emph{2-4-2} building block, so that the bottom row of the polygon has
the same length as the top row of the building block. This
construction is depicted in Figure~\ref{fig 242 recomp}.

We will translate this construction into a recurrence satisfied by
the \emph{2-4-2} polygon generating function by using Hadamard
products. We note that a similar construction (but for different
lattice objects) appears in \cite{klarner-can, Klarner1968} but is
phrased in terms of constant term integrals.

\begin{figure}[ht!]
  \centering
  \includegraphics[scale=0.5]{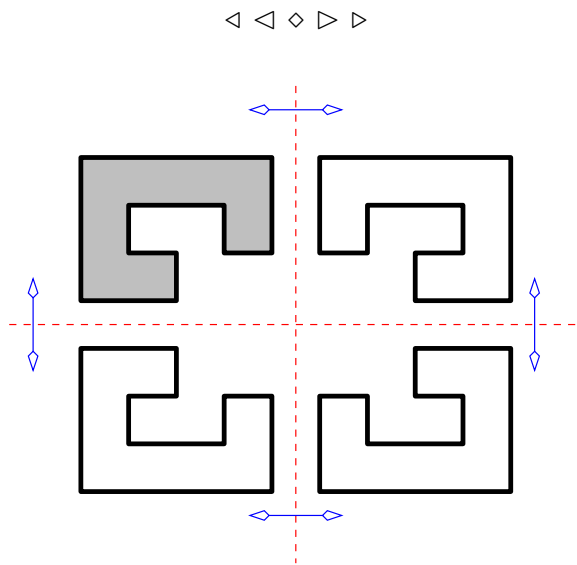}
  \caption{The set of building blocks has a 4-fold symmetry as shown. It suffices to find all the building blocks in one orientation and then obtain the others by reflections.}
  \label{fig 242 orient}
\end{figure}

Let us start with the generating function of the building blocks:
\begin{lemma}
  Let $T(t,s;x,y)$ be the generating function of \emph{2-4-2} polygon
  building blocks, where $t$ and $s$ are conjugate to the length of
  top and bottom rows (respectively). Then $T$ may be expressed as 
  \begin{equation}
    T(t,s;x,y) = 2\left(\hat{T}(t,s;x,y) + \hat{T}(s,t;x,y) \right),
  \end{equation}
  where the generating function $\hat{T}(t,s;x,y)$ is given by
\begin{eqnarray}
  \label{eqn BB_fixed}
  \hat{T}(t,s;x,y) & = & y^4 \Bigl( A(s,t;x)\cdot \seq{stx}\seq{tx}^2 \cdot B(s,t;x) \Bigr. \nn 
  & & + A(s,t;x)\cdot\seq{stx}\seq{stx^2}\seq{tx}^2 \cdot B(s,t;x) \nn 
  & & + A(s,t;x)\cdot \seq{stx}\seq{tx}^3 \cdot B(s,t;x) \nn 
  & & + C(s,t;x)\cdot \seq{sx}\seq{tx}^3 \cdot B(s,t;x) \nn 
  & & + \Bigl. C(s,t;x)\cdot \seq{sx}\seq{x}\seq{tx}^3 \cdot B(s,t;x) \Bigr).
\end{eqnarray}
We have used $\seq{f}$ as shorthand for $\frac{f}{1-f}$, and the
generating functions $A$, $B$ and $C$ are:
\begin{eqnarray}
  \label{eqn frills}
  A(s,t;x) & = & 1 + \seq{x} + 2 \seq{sx} + 2\seq{tx} +  \seq{sx}\seq{tx} + \nonumber \\ 
  & &   \seq{sx}^2 + \seq{sx}\seq{x}  + \seq{tx}^2 + \seq{tx}\seq{x} \\ 
  B(s,t;x) & = & 1 + \seq{tx} +
  \seq{x} \\ C(s,t;x) & = & 1 + \seq{sx} + \seq{x}.
\end{eqnarray}
\end{lemma}
\proof  Figure~\ref{fig 242 orient} shows the four possible
orientations of a building block. Figures~\ref{fig 242
  bb}~and~\ref{fig 242 frills} show how to construct the generating
function $\hat{T}$ of building blocks in one orientation. To obtain
all building blocks we must reflect the blocks counted by $\hat{T}$
about both horizontal and vertical lines (as shown in Figure~\ref{fig
  242 orient}). Reflecting about a vertical line multiplies $\hat{T}$ by
$2$. Reflecting about a horizontal line interchanges the roles of $s$
and $t$. This proves the first equation.

\begin{figure}[ht!]
  \centering
  \includegraphics[scale=0.5]{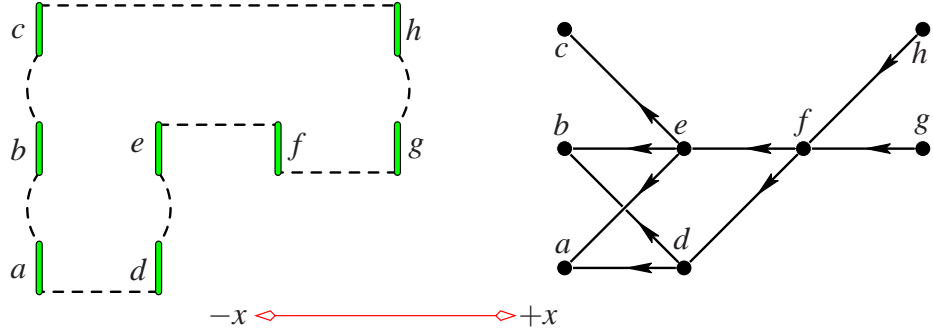}
  \caption{The vertical bonds of a \emph{2-4-2} polygon building
    block. The $x$-ordinate of these bonds are denoted $a,b,\dots,h$
    as shown. The Hasse diagram showing the constraints on the values
    $a, b, \dots, h$ is given on the right; an arrow from $v_i$ to $v_j$
    implies that $v_i>v_j$.}
  \label{fig 242 hasse}
\end{figure}

We now find $\hat{T}$ by finding the \emph{section-minimal} building
blocks in one orientation (that of the top-left polygon in
Figure~\ref{fig 242 orient}).  All such polygons contain $8$ vertical
bonds, let $a,\dots,h \in \mathbb{Z}$ denote the $x$-ordinate of these
bonds.  Figure~\ref{fig 242 hasse} shows the Hasse diagram that these
numbers must satisfy:
\begin{displaymath}
  \begin{array}{llr}
    a,b < d \qquad & a,b,c < e \\
    d,e < f &  f < g,h.
\end{array}
\end{displaymath}
Without loss of generality we set $a=0$ (to enforce translational
invariance). 

\begin{figure}[ht!]
  \centering 

  \psfrag{t}{$t$} \psfrag{s}{$s$} \psfrag{A}{$A$}
  \psfrag{B}{$B$} \psfrag{C}{$C$}

  \psfrag{1}{1}\psfrag{2}{2}\psfrag{3}{3}
  \psfrag{4}{4}\psfrag{5}{5}

  \includegraphics[scale=0.5]{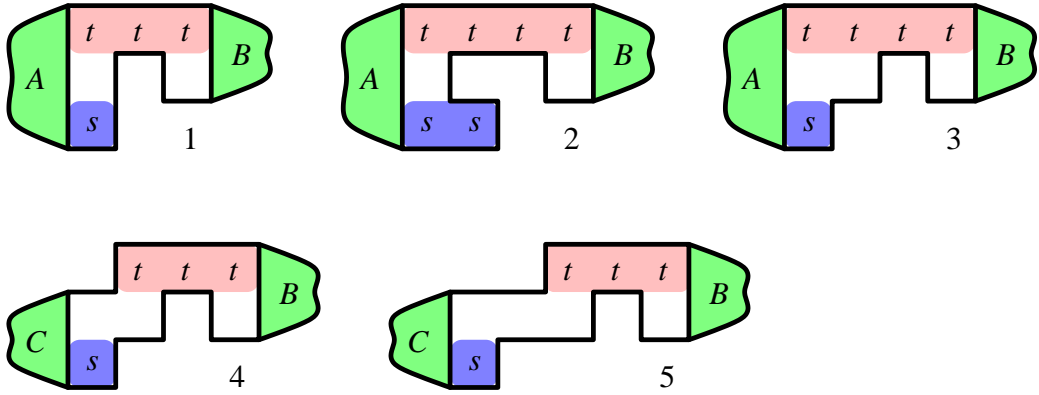}
  \caption{The section-minimal building blocks of \emph{2-4-2} polygons. The ``\emph{frills}'', denoted $A$, $B$ and $C$ are given in Figure~\ref{fig 242 frills}. }
  \label{fig 242 bb}
\end{figure}

\begin{figure}[ht!]
  \centering \psfrag{t}{$t$} \psfrag{s}{$s$} \psfrag{A}{$A$}
  \psfrag{B}{$B$} \psfrag{C}{$C$}
  \includegraphics[scale=0.5]{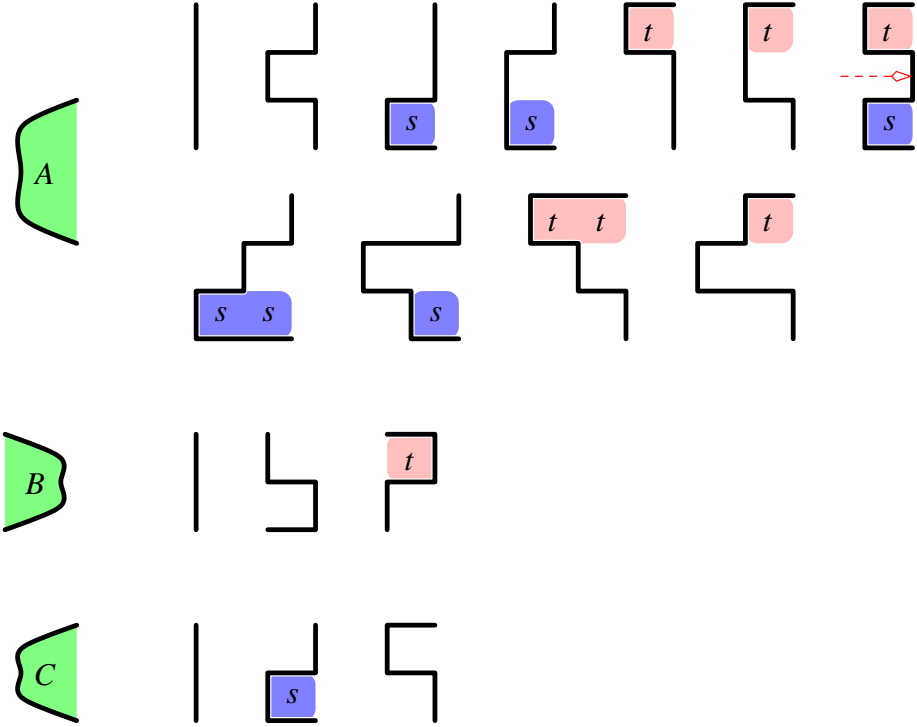}
  \caption{The ``\emph{frills}'' of the building blocks in Figure~\ref{fig 242 bb}.}
  \label{fig 242 frills}
\end{figure}

Consider a section-minimal building block and determine the values of
$b, \dots, h$. We can decompose the building block depending on these
values:
\begin{itemize}
\item Find which of $g$ and $h$ is minimal and cut the polygon along a vertical
  line running through that vertical bond. This separates the polygon
  into 2 parts; the part to the right is a \emph{B-frill} (see
  Figure~\ref{fig 242 frills}) --- there are 3 possible
  \emph{B-frills} depending on whether $g=h$, $g<h$ or $g>h$.
  
\item If $c<d$ then the building block must be of the form of polygon
  $1$, $2$ or $3$ in Figure~\ref{fig 242 bb}. Determine which is the
  greatest of $a,b$ and $c$ and cut the polygon along the vertical line
  running through that vertical bond. This separates the polygon into
  2-parts; the part to the right is an \emph{A-frill} (see
  Figure~\ref{fig 242 frills}) --- there are
  11 possible \emph{A-frills} depending on the relative magnitudes
  of $a$, $b$, and $c$.

\item If $c \geq d$ then the building block must be of the form of polygon
  $4$ or $5$ in Figure~\ref{fig 242 bb}. Find which of $a$ and $b$ is
  greater and cut along the vertical line running through that
  vertical bond. This separates the polygon into 2 parts; the part to
  the right is a \emph{C-frill} (see Figure~\ref{fig 242 frills}) ---
  there are 3 possible \emph{C-frills} depending on the whether
  $a=b$, $a<b$ or $a>b$.
\end{itemize}
Using this decomposition we see that every section minimal polygon is
given by one of the 5 polygons given in Figure~\ref{fig 242 bb}
together with 2 of the \emph{frills} from Figure~\ref{fig 242 frills}. The
above equation for $\hat{T}(t,s;x,y)$ follows. \qed

We note that one could find $\hat{T}$ using the theory of
$P$-partitions \cite{stane}, and we used it to check the result. 

\spacebreak




We now define the (restricted) Hadamard product and show how it
relates to the construction of \emph{2-4-2} polygons.

\begin{defn}
  Let $f(t) = \sum_{t\geq0} f_n t^n$ and $g(t) = \sum_{t\geq0} g_n
  t^n$ be two power series in $t$. We define the (restricted) Hadamard
  product $f(t) \hp{t} g(t)$ to be
  \begin{displaymath}
    f(t) \hp{t} g(t) = \sum_{n \geq 0} f_n g_n.
  \end{displaymath}
\end{defn}
We note that if  $f(t)$ and $g(t)$ are two power series with real
coefficients such that 
\begin{displaymath}
  \lim_{n \to \infty} \left| f_n g_n \right|^{1/n} < 1,    
\end{displaymath}
then the Hadamard product $f(t) \hp{t} g(t)$ will exist. For example
$(1-2t)^{-1} \hp{t} (1-3t)^{-1}$ does not exist, while $(1-2t)^{-1}
\hp{t} (1-t/3)^{-1}$ does exist and is equal to 3.

Below we consider Hadamard products of power series in $t$ whose
coefficients are themselves power series in two variables, $x$ and
$s$. These products are of the form
\begin{equation}
  f(t;x) \hp{t} T(t,s,x) = \sum_{n \geq 0} f_n(x) T_n(s,x).
\end{equation}
Since the summands are the generating functions of certain polygons
(see below) it follows that $f_n(x) T_n(x) = O(s x^n)$ and so the sum
converges.

\begin{lemma}
  Let $f(s;x,y)$ be the generating function of \emph{2-4-2}~polygons,
  where $s$ is conjugate to the length of bottom row of the
  polygon. This generating function satisfies the following equation
  \begin{displaymath}
    f(s;x,y) = \frac{ysx}{1-sx} + 
    f(t;x,y) \hp{t} \left(\frac{1}{y} T(t/x,s;x,y) \right),
  \end{displaymath}
  where $T(t,s;x,y)$ is the generating function of the \emph{2-4-2}
  building blocks.
\end{lemma}
\proof Let us write $f(s;x,y) = \sum_{n \geq 1} f_n(x,y) s^n$ and
$T(t,s;x,y) = \sum_{n \geq 1} T_n(s;x,y)$, where $f_n(x,y)$ is the
generating function of \emph{2-4-2} polygons whose \emph{bottom} row
has length $n$, and $T_n(s;x,y)$ is the generating function of
\emph{2-4-2} building blocks, whose \emph{top} row has length $n$. The
above recurrence becomes:
\begin{displaymath}
  f(s;x,y) = \frac{ysx}{1-sx} 
  + \sum_{n \geq 1}  f_n(x,y) T_n(s;x,y) / (yx^n).
\end{displaymath}
This follows because \emph{2-4-2} polygon is either a rectangle of
unit height (counted by $\frac{ysx}{1-sx}$) or may be constructed by
combining a \emph{2-4-2} polygon, whose last row is of length $n$
(counted by $f_n(x,y)$) with a \emph{2-4-2} polygon whose top row is
of length $n$ (counted by $T_n(s;x,y)$). To explain the factor of
$1/(yx^n)$ see Figure~\ref{fig 242 recomp}; when the building block is
joined to the polygon (centre) and the duplicated row is
``\emph{squashed}'' (right), the total vertical half-perimeter is
reduced by $1$ (two vertical bonds are removed) and the total
horizontal half-perimeter is reduced by the length of the join (two
horizontal bonds are removed for each cell in the join). Hence if the
join is of length $n$, the perimeter weight needs to be reduced by a factor
of $(yx^n)$. \qed

\spacebreak

While in general Hadamard products are difficult to evaluate, if one
of the functions is rational then the result is quite simple.  This
fact allows us to translate the above Hadamard-recurrence into a
functional equation.
\begin{lemma}
  \label{lem hadamard}
  Let $f(t) = \sum_{t \geq 0} f_n t^n$ be a power series. The
  following (restricted) Hadamard products are easily evaluated:
  \begin{eqnarray}
    f(t) \hp{t} \left( \frac{1}{1-\alpha t} \right) & = & f(\alpha) \\ 
    f(t) \hp{t} \left( \frac{k! t^k}{(1-\alpha t)^{k+1}} \right) & = & 
    \left.\DDiff{f}{t}{k}\right|_{t =\alpha}.
  \end{eqnarray}
  We also note that the Hadamard product is linear:
  \begin{equation}
    f(t) \hp{t} \big( g(t) + h(t) \big) = f(t) \hp{t} g(t) + f(t) \hp{t} h(t).
  \end{equation}
\end{lemma}
\proof The second equation follows from the first by differentiating
with respect to $\alpha$. The first equation follows because
\begin{displaymath}
  f(t) \hp{t} \frac{1}{1-\alpha t}  =
  \left( \sum_{n \geq 0} f_n t^n \right) \hp{t} \left(\sum_{n \geq 0} \alpha^n t^n \right)
  = \sum_{n \geq 0} f_n \alpha^n = f(\alpha).
\end{displaymath}
The linearity follows directly from the definition.
\qed

\spacebreak

In order to apply the above lemma, we need to rewrite $T(t/x,s;x,y)/y$ in
(a non-standard) partial fraction form:
\begin{equation}
  T(t/x,s;x,y)/y = y^3 \left[ c_0 \cdot t^0 + \sum_{k=0}^{5} c_{k+1}
    \frac{ k! t^k}{(1-t)^{k+1}} 
  + c_7 \frac{1}{1-st} + c_8 \frac{1}{1-stx} \right],
\end{equation}
where the $c_i$ are rational functions of $s$ and $x$. We note
that when $s=1$ some singularities of $T$ coalesce and we rewrite $T$
as:
\begin{equation}
  T(t/x,1;x,y)/y = y^3 \left[ \hat{c}_0 \cdot t^0 
    + \sum_{k=0}^{6} \hat{c}_{k+1} \frac{ k! t^k}{(1-t)^{k+1}} 
  + \hat{c}_8 \frac{1}{1-tx} \right],
\end{equation}
where the $\hat{c}_i$ are rational functions of $x$. The Hadamard
product $f(t;x,y) \hp{t} T(t/x,s;x,y)/y$ is then:
\begin{equation}
  \label{eqn 242 s=s}
  f(t;x,y) \hp{t} T(t/x,s;x,y)/y = y^3 \left[ 
    \sum_{k=0}^{5} c_{k+1} \Ddiff{f}{t}{k}(1;x,y) + c_7 f(s;x,y) + c_8 f(sx;x,y) \right],
\end{equation}
where we have made use of the fact that $[t^0] f(t;x,y) = 0$ (there are no rows of zero length).

We do not state in full the coefficients, $c_i$, since they are very
large and, with the exception of $c_8$, not particularly relevant to
the following analysis. We will just state the denominators of all the
coefficients, as well as the coefficient $c_8$ in full.  If we write
the denominator of $c_i$ as $d_i$:
\begin{eqnarray}
  d_0 = (1-x)^3 (1-sx)^6 (1-s)^6  &&  
  d_1  =  (1-x)^3 (1-sx)^5 (1-s)^5 \nonumber \\ 
  d_2 = (1-x)^3 (1-sx)^3 (1-s)^4 &&
  d_3  = (1-x)^3 (1-sx)^3 (1-s)^3 \nonumber \\ 
  d_4 =  (1-x)^2 (1-sx) (1-s)^2 &&
  d_5 = (1-x) (1-sx) (1-s) \nonumber \\ 
  d_6  = (1-s) &&  
  d_7 = (1-sx)^6 (1-s)^6 \nonumber \\ 
  c_8 = - \frac{2sx^2(s^2x^2+sx-s+1)}{(1-sx)^4 (1-x)^2}. &&
\end{eqnarray}

When $s=1$ the coalescing poles change equation~\Ref{eqn 242 s=s} to:
\begin{equation}
  \label{eqn 242 s=1}
  f(t;x,y) \hp{t} T(t/x,1;x,y)/y = y^3\left[ 
    \sum_{k=0}^{6} \hat{c}_{k+1} \Ddiff{f}{t}{k}(1;x,y) + \hat{c}_8 f(x;x,y) \right]
\end{equation}
The coefficients, $\hat{c}_i$, become somewhat simpler and can be stated here in full:
\begin{eqnarray}
  \hat{c}_0 = -2\frac{x^3(1+x)(2x^2+1)}{(1-x)^6} &&
  \hat{c}_1 = 4\frac{(1+x)(x^2+1)x^3}{(1-x)^6} \nonumber \\
  \hat{c}_2  =  2\frac{x^2(1+x)(2x^2+x+1)}{(1-x)^5} &&
  \hat{c}_3  =  \frac{x^2(1+x)(2x+1)}{(1-x)^4} \nonumber \\
  \hat{c}_4  =  \frac{1}{3} \frac{(1+x)(x^2+x+1)}{(1-x)^3} &&
  \hat{c}_5  =  \frac{1}{12}\frac{(x^2+2x+3)}{(1-x)^2} \nonumber \\ 
  \hat{c}_6  =  \frac{1}{60}\frac{(x+3)}{(1-x)} &&
  \hat{c}_7  =  \frac{1}{360} \nonumber \\ 
  \hat{c}_8  =  -2\frac{x^3(1+x)}{(1-x)^6} & = & c_8|_{s=1}.
\end{eqnarray}

\begin{lemma}
  \label{lem 242 feqn}
  Let $f(s;x,y)$ be the generating function for \emph{2-4-2}~polygons
  enumerated by bottom row-width, half-horizontal perimeter and
  half-vertical perimeter ($s, x$ and $y$ respectively).  $f(s;x,y)$
  satisfies the following functional equations:
  \begin{eqnarray}
    \label{eqn 242 func eqns}
    f(s;x,y) &=& \frac{sxy}{1-sx} + y^3\left[ \sum_{k=0}^{5} c_{k+1} \Ddiff{f}{s}{k}(1;x,y) 
      + c_7 f(s;x,y) + c_8 f(sx;x,y) \right] \qquad \\ 
    f(1;x,y) &=& \frac{xy}{1-x} + y^3\left[ \sum_{k=0}^{6} \hat{c}_{k+1} \Ddiff{f}{s}{k}(1;x,y) 
      + \hat{c}_8 f(x;x,y) \right],
  \end{eqnarray}
  with $c_i$ and $\hat{c}_i$ given above.
  
  We rewrite the generating function as $f(s;x,y) = \sum_{n\geq 1} f_n(s;x) y^{3n-2}$, where
  the coefficient $f_n(s;x)$ is the generating function for ${\cal
    P}^{242}_{3n-2}$. This allows the above functional equations to be
  transformed into recurrences:
  \begin{eqnarray}
    f_1(s;x) & = & \frac{sx}{1-sx} \\ 
    f_{n+1}(s;x) & = & \sum_{k=0}^{5} c_{k+1} \Ddiff{f_n}{s}{k}(1;x) 
    + c_7 f_n(s;x) + c_8 f_n(sx;x) \qquad  s \neq 1 \\ 
    f_{n+1}(1;x) & = & \sum_{k=0}^{6} \hat{c}_{k+1} \Ddiff{f_n}{s}{k}(1;x) + \hat{c}_8 f_n(x;x).
  \end{eqnarray}
\end{lemma}
\proof Apply Lemma~\ref{lem hadamard} to the partial fraction form of the transition
function for general $s$, and when $s=1$.  \qed

\subsection{Analysing the functional equation}
By Lemma~\ref{lem 242 pole}, we are able to prove Theorem~\ref{thm
  first ksec exp} by showing that $f_n(1;x)$ is singular at the zeros
of $\Psi_n(x)$. We do this by induction using the recurrences in
Lemma~\ref{lem 242 feqn}.

Before we can do this we need to prove the following lemma about the
zeros (and hence factors) of one of the coefficients in the recurrence:
\begin{lemma}
  \label{lem c8 zero}
  Consider the coefficient $c_8(s;x)$ defined above. When $s=x^k$,
  $c_8(x^k,x)$ has a single zero on the unit circle at $x=-1$ when $k$
  is even. When $k$ is odd $c_8(x^k,x)$ has no zeros on the unit circle.

\end{lemma}
\proof When $s=x^k$,  $c_8(x^k,x)$ is
\begin{displaymath}
  c_8(x^k,x) = \frac{ 2 x^{k+2}( k^{2k+2} + x^{k+1} - x^k + 1) }
  { (1-x^{k+1})^4(1-x)^2 }.
\end{displaymath}
Let $\xi$ be a zero of $c_8(x^k,x)$ that lies on the unit circle;
$\xi$ must be a solution of the polynomial $p_k(x) = x^{2k+2} + x^{k+1}
- x^k + 1 = 0$. Hence:
\begin{eqnarray*}
  \xi^k - \xi^{k+1} & = & \xi^{2k+2} + 1 \qquad \qquad \mbox{ divide by $\xi^{k+1}$}\\
  1/\xi - 1 & = & \xi^{k+1} + \xi^{-k-1}.
\end{eqnarray*}
Since $\xi$ lies on the unit circle we may write $\xi = e^{i \theta}$:
\begin{eqnarray*}
  e^{- i \theta}  - 1  & = &  e^{i (k+1) \theta} + e^{ - i (k+1) \theta} \\
  & = & 2 \cos( (k+1) \theta ).
\end{eqnarray*}
Since the right hand-side of the above expression is real the
left-hand side must also be real. Therefore $\theta = 0, \pi$ and $\xi
= \pm 1$. If $\xi = 1$ then $p_k(\xi) = 2$. On the other hand, if $\xi
= -1$ then $p_k(\xi) = 4$ if $k$ is odd and is zero of $k$ is even.

Since the denominator of $c_8(x^k,x)$ is not zero when $k$ is even and
$x=-1$ the result follows. One can verify that there are not multiple
zeros at $x=-1$ by examining the derivative of the numerator.

\qed

\spacebreak

\noindent \textbf{Proof of Theorem~\ref{thm first ksec exp}}:\\
This proof for SAPs was first given in \cite{ADR_thesis}. A similar
(but cleaner) argument for a different class of polygons appears in
\cite{MBM_ADR_02}. We follow the latter.

Consider the recurrence given in Lemma~\ref{lem 242 feqn}.  This
recurrence shows that $f_n(s;x)$ may be written as a rational function
of $s$ and $x$. Since $f_n(1;x)$ is a well defined (and rational)
function, the denominator of $f_n(s;x)$ does not contain any factors
of $(1-s)$.

Let $\mathbb{C}_n(s;x)$ be the set of polynomials of the form
\begin{equation}
  \prod_{k=1}^{n} \Psi_k(x)^{a_k} (1-sx^k)^{b_k},
\end{equation}
where $a_k$ and $b_k$ are non-negative integers. We define
$\mathbb{C}_n(x) = \mathbb{C}_n(0;x)$ (polynomials which are products
of cyclotomic polynomials). We first prove (by induction on $n$) that
$f_n$ may be written as
\begin{equation}
  \label{eqn fn denom form}
  f_n(s;x) = \frac{N_n(s;x)}{ (1-sx^n) D_n (s;x) },
\end{equation}
where $N_n(s;x)$ and $D_n(s;x)$ are polynomials in $s$ and $x$ with
the restriction that $D_n(s;x) \in \mathbb{C}_{n-1}(s;x)$. Then we
consider what happens at $s=1$ and $x$ is set to a zero of $\Psi_k$.

For $n=1$, equation~\Ref{eqn fn denom form} is true, since $f_1(s;x) =
\frac{sx}{1-sx}$. Now assume equation~\Ref{eqn fn denom form} is true
up to $n$ and apply the recurrence. The only term that may introduce a
new zero into the denominator is $c_8(s;x) f_n(sx;x)$. By
assumption $f_n(s;x) = \frac{N_n(sx;x)}{(1-sx^{n+1}) D_n(sx;x)}$, and
$D_n(sx;x) \in C_{n}(s;x)$. Hence equation~\Ref{eqn fn denom form} is
true for $n+1$, and so is also true for all $n \geq 1$.

\spacebreak

Let $\xi$ be a zero of $\Psi_k(x)$. We wish to prove that $f_n(1;x)$
is singular at $x=\xi$ and we do so by proving that for $k=1, \dots,
n$, the generating function $f_k(x^{n-k};x)$ is singular at $x=\xi$,
and then setting $k=n$. We proceed by induction on $k$ for fixed $n$.

If we set $k=1$, then we see that $f_1(x^{n-1};x) =
\frac{x^n}{1-x^n}$, and so the result is true. Now let $k \geq 2$ and
assume that the result is true for $k-1$, \emph{ie}
$f_{k-1}(x^{n-k+1};x)$ is singular at $x=\xi$. The recurrence relation
and equation~\Ref{eqn fn denom form} together imply
\begin{equation}
  f_k(s; x) = \frac{N(s;x)}{D(s;x)}  
  + c_8(s;x) f_{k-1}(sx;x),
\end{equation}
where $N$ and $D$ are polynomials in $s$ and $x$ and $D(s;x) \in
\mathbb{C}_{k-1}(s;x)$. Setting $s=x^{n-k}$ yields
\begin{equation}
  \label{eqn sing rec}
  f_k(x^{n-k};x) = \frac{N(x^{n-k};x)}{D(x^{n-k};x)} 
  + c_8(x^{n-k};x) f_{k-1}(x^{n-k+1};x),
\end{equation}
and we note that $D(x^{n-k};x) \in \mathbb{C}_{n-1}(x)$. In the case
$k=n$ the above equation is still true, since $\hat{c}_8 = c_8|_{s=1}$.

Equation~\Ref{eqn sing rec} shows that $f_k(x^{n-k})$ is singular at
$x = \xi$ only if $c_8(x^{n-k};x) f_{k-1}(x^{n-k+1};x)$ is singular at
$x = \xi$. This is true (by assumption) unless $c_8(x^{n-k};x) = 0$ at
$x=\xi$. By Lemma~\ref{lem c8 zero}, $c_8(x^{n-k};x)$ is non-zero at
$x=\xi$, except when $n=k=2$.

In the case $n=k=2$ this proof breaks down, and indeed we see that
$H_4(x)$ is not singular at $x=-1$. Excluding this case,
$f_k(x^{n-k};x)$ is singular at $x = \xi$ and so $f_n(1;x)$ is also
singular at $x = \xi$. By Lemma~\ref{lem 242 pole}, $H_{3k-2}(x)$ is
singular at $x = \xi$. \qed

\spacebreak
We can now prove that the self-avoiding polygon anisotropic generating
function is not a D-finite function:
\begin{cor}
  \label{cor sap ndf}
  Let $S_n$ be the set of singularities of the coefficient $H_n(x)$.
  The set $S = \bigcup_{n\geq1} S_n$ is dense on the unit circle $|x|=1$.
  Consequently the self-avoiding polygon anisotropic half-perimeter
  generating function is not a D-finite function of $y$.
\end{cor}
\proof For any $q \in \mathbb{Q}$, there exists $k$, such that
$\Psi_k( e^{2\pi i q}) = 0$. By Theorem~\ref{thm first ksec exp},
$H_{3k-2}(x)$ is singular at $x = e^{2\pi i q}$, excepting
$x=-1$. Hence the set $S$ is dense on the unit circle,
$|x|=1$. Consequently $S$ has an infinite number of accumulation
points and so $G(x,y) = \sum H_n(x) y^n$ is not a D-finite power
series in $y$. 

\qed

We can easily extend this result to self-avoiding polygons on
hypercubic lattices.

\begin{cor}
  \label{cor hypercubic}
  Let $\mathcal{G}_d$ be the set of self-avoiding polygons on the
  $d$-dimensional hypercubic lattice and let $G_d$ be the anisotropic
  generating function
  \begin{displaymath}
    G_d( x_1, \dots, x_{d-1}, y) = 
    \sum_{P \in \mathcal{G}_d} y^{|P|_d}\prod_{i=1}^{d-1} x_i^{|P|_i},
  \end{displaymath}
  where $|P|_i$ is half the number of bonds in parallel to the unit
  vector $\vec{e}_i$. \emph{Ie} when $d=2$ we recover the square
  lattice anisotropic generating function. If $d=1$ then the
  generating function is equal to zero and otherwise is a
  non-D-finite power series in $y$.
\end{cor}
\proof When $d=1$ then there are no self-avoiding polygons and so the
generating function is trivially zero. Now consider $d \geq 2$. It is
a standard result in the theory of D-finite power series that any well
defined specialisation of a D-finite power series is itself D-finite
\cite{Lipshitz1989}. Setting $x_2 = \dots = x_{d-1} = 0$ in the
generating function $G_d(\mathbf{x},y)$ recovers the square lattice
generating function $P(x,y)$. Hence if $G_d(\mathbf{x}, y)$ were a
D-finite power series in $y$ then it follows that $P(x,y)$ would also
be D-finite. This contradicts Corollary~\ref{cor sap ndf}, and so the
result follows.

\qed

\section{Discussion}
We have shown that the anisotropic generating function of
self-avoiding polygons on the square lattice, $P(x,y)$, is not a
D-finite function of $y$. This result was then extended to prove that
the anisotropic generating function of self-avoiding polygons on any
hypercubic lattice is either trivial (in one dimension) or a
non-D-finite function (in dimensions 2 and higher).

There exists a number of non-D-finiteness results for generating
functions of other models, such as bargraphs enumerated by their
site-perimeter \cite{MBM_ADR_02}, a number of lattice animal models
related to heaps of dimers \cite{MBM_AR_CH} and certain types of
matchings \cite{Klazar}; these results rely upon a knowledge of the
generating function --- either in closed form or via some sort of
recurrence. The result for self-avoiding polygons is, as far as we are
aware, the first result on the D-finiteness of a completely unsolved
model!

Unfortunately we are not able to use this result to obtain information
about the nature of the isotropic generating function $P(x,x)$; it is
all too easy to construct a two-variable function that is not D-finite
, that reduces to a single variable D-finite function. For example
\begin{equation}
  F(x,y) = \sum_{n \geq 1} \frac{y^n}{(1-x^n)(1-x^{n+1})}.
\end{equation}
is not a D-finite function of $y$ by Theorem~\ref{thm Dfinite
  poles}. Setting $y=x$ gives a rational, and hence D-finite, function
of $x$:
\begin{eqnarray*}
  F(x,x) & = & \sum_{n \geq 1} \frac{x^n}{(1-x^n)(1-x^{n+1})} \\
  & = & \frac{1}{1-x} \sum_{n \geq 1} 
  \left( \frac{x^{n}}{1-x^{n}} - \frac{x^{n+1}}{1-x^{n+1}} \right) \\
  & = & \frac{x}{(1-x)^2}.
\end{eqnarray*}

On the other hand, the anisotropisation of \emph{solvable} lattice
models does not alter the nature of the generating function --- rather
it simply moves singularities around in the complex plane.
Unfortunately we are unable to rigorously determine how far this
phenomenon extends since we know very little about the nature of the
generating functions of unsolved models.

That the self-avoiding polygon anisotropic half-perimeter generating
function is not D-finite (Corollary~\ref{cor sap ndf}) demonstrates
the stark difference between the bond-animal models we have been able
to solve and those we wish to solve. Solved bond-animal models (with
the exception of spiral walks \cite{blote84}) all have D-finite
anisotropic generating functions. More general (and unsolved) models,
such as bond animals and self-avoiding walks, are believed to exhibit
the same dense pole structure \cite{Guttmann2000, Guttmann1996} as
self-avoiding polygons and therefore are thought to be non-D-finite.

Two papers are in preparation to extend these results to directed bond
animals, bond trees and general bond animals. We are also
investigating the possibility of applying these techniques to
site-animals and other combinatorial objects.


\section*{Acknowledgements}
I would like to thank a number of people for their assistance during
this work.
\begin{itemize}
\item A.\ L.\ Owczarek, A.\ J.\ Guttmann and M.\ Bousquet-M\'elou for
  supervision and advice during my PhD and beyond.
\item I.\ Jensen for his anisotropic series data.
\item R.\ Brak, M.\ Zabrocki, E.\ J.\ Janse van Rensburg and N.\ Madras for
  discussions and help with the manuscript.
\item The people at LaBRI for their hospitality during my stay at the Universit\'e Bordeaux 1.
\end{itemize}
This work was partially funded by the Australian Research Council.


\bibliographystyle{plain}
\bibliography{dfinite.bib}

\end{document}